\documentclass[11pt]{amsart} 
\usepackage[all]{xy}
\usepackage{shuffle,tikz,stmaryrd,stackrel,mathtools,amsmath,amssymb,mathrsfs}

\newcommand{\C}{\mathbb C}
\newcommand{\N}{\mathbb N}
\newcommand{\Z}{\mathbb Z}
\newcommand{\ssl}{\mathfrak{sl}}
\newcommand{\fg}{\mathfrak{g}}

\newcommand{\ft}{\mathfrak t}
\newcommand{\fn}{\mathfrak n}
\newcommand{\eps}{\varepsilon}
\newcommand\SL{{\mathbf{SL}}}
\newcommand\SO{{\mathbf{SO}}}
\newcommand{\dimvec}{\underrightarrow{\dim}\,}
\newcommand{\vi}{{\underline i}}
\newcommand{\vj}{{\underline j}}
\newcommand{\vk}{{\underline k}}
\newcommand{\tR}{{\ft_{\mathbb R}^*}}
\newcommand{\R}{\mathbb R}

\newcommand{\shin}{\shuffle}

\newcommand{\barD}{\overline{D}}
\newcommand\mmod{\text{\upshape-}\mathrm{mod}}
\newcommand{\ftreg}{{\ft^{\mathrm{reg}}}}

\DeclareMathOperator\im{im}
\DeclareMathOperator\wt{wt}
\DeclareMathOperator\Hom{Hom}
\DeclareMathOperator\Conv{Conv}
\DeclareMathOperator\Gr{Gr}

\DeclareMathOperator\Pol{Pol}
\DeclareMathOperator\Irr{Irr}

\DeclareMathOperator\Ad{Ad}
\DeclareMathOperator\Rep{Rep}
\DeclareMathOperator\Dist{Dist}
\DeclareMathOperator\ttop{top}

\theoremstyle{plain}
\newtheorem{Theorem}{Theorem}[section]
\newtheorem{Lemma}[Theorem]{Lemma}
\newtheorem{Corollary}[Theorem]{Corollary}
\newtheorem{Problem}[Theorem]{Problem}
\newtheorem{Question}[Theorem]{Question}
\newtheorem{Proposition}[Theorem]{Proposition}
\newtheorem{Conjecture}[Theorem]{Conjecture}
\theoremstyle{definition}
\newtheorem{Example}[Theorem]{Example}
\newtheorem{Remark}[Theorem]{Remark}

\numberwithin{equation}{section}

\begin{document}


\title[Perfect bases in representation theory]{Perfect bases in representation theory: three mountains and their springs}

\author{Joel Kamnitzer}


\address{Department of Mathematics, University of Toronto, Toronto ON, Canada} \email{jkamnitz@math.toronto.edu}


\begin{abstract}
In order to give a combinatorial descriptions of tensor product multiplicites for semisimple groups, it is useful to find bases for representations which are compatible with the actions of Chevalley generators of the Lie algebra.  There are three known examples of such bases, each of which flows from geometric or algebraic mountain.  Remarkably, each mountain gives the same combinatorial shadow: the crystal $B(\infty) $ and the Mirkovi\'c--Vilonen polytopes.  In order to distinguish between the three bases, we introduce measures supported on these polytopes.  We also report on the interaction of these bases with the cluster structure on the coordinate ring of the maximal unipotent subgroup.
\end{abstract}

\maketitle
 

\section{Representations and their bases}

\subsection{Semisimple Lie algebras and their representations}

Let $ G $ be a complex semisimple group. The representation theory of $ G $ is very well understood.  The irreducible $ G$-rep\-resen\-tations are labelled by dominant weights, and every representation is a direct sum of these irreducible representations.  For $ \lambda \in P_+$, the irreducible representation $ V(\lambda) $ admits a decomposition into eigenspaces $V(\lambda)_\mu $ for the action of $ T $.  These eigenspaces are called weight spaces and their dimensions are called \textbf{weight multiplicities}.

The tensor product of two irreducible representations decomposes into a direct sum of irreducible representations with \textbf{tensor product multiplicities} $c^\nu_{\lambda \mu} $.
$$
V(\lambda) \otimes V(\mu) \cong \bigoplus_{\nu \in P_+} V(\nu)^{\oplus c^\nu_{\lambda \mu}}
	$$
	
	\begin{Problem} \label{prob:mult}
		Determine combinatorial formulae for weight multiplicities and tensor product multiplicities. 
	\end{Problem}

Weight and tensor product multiplicities are closely related by the following construction.  Let $ C^\nu_{\lambda \mu} = \Hom(V(\nu), V(\lambda) \otimes V(\mu))$, a vector space whose dimension is $ c^\nu_{\lambda \mu}$.  
\begin{Proposition} \label{pr:Zhelo}
	There is an injective map $ C^\nu_{\lambda \mu} \rightarrow V(\lambda)_{\nu - \mu} $ with image $ \displaystyle{\bigcap_{i \in I}} \ker \, e_i^{\alpha_i^\vee(\mu) + 1}$.
\end{Proposition}
Here we use the Chevalley presentation of $ \fg $, with generators $ e_i, f_i, \alpha_i^\vee $, for $ i \in I $.

	\subsection{Good and perfect bases}
	Problem \ref{prob:mult} was first solved by Littelmann \cite{L} and Berenstein--Zelevinsky \cite{BZ}, following an approach first proposed by Gel$'$fand--Zelevinsky \cite{GZ}.  They suggested finding weight bases for each $ V(\lambda) $ which restrict to bases of tensor product multiplicity spaces. 
	
	Let $ V $ be a $G$-representation.  A \textbf{weight basis} for $ V $ is a basis consisting of weight vectors.  A weight basis $ B $ for $ V(\lambda) $ is called \textbf{good}, if for each $ i \in I$, it is compatible with the filtration of $ V(\lambda) $ given by the kernels of powers of $ e_i $.  From Proposition \ref{pr:Zhelo}, it follows that a good basis restricts to a basis of each tensor product multiplicity space.
	
	
	A slight strengthening of the notion of good basis was proposed by Berenstein--Kazhdan \cite{BK}.  One might imagine that we could find a basis for a representation such that each $ e_i $ takes each basis vector to another basis vector (or 0).  However, this is not always possible (see Example \ref{eg:sl3adjoint}).  So instead we will demand that each $ e_i $ permutes the basis up to lower order terms.
	
	To formulate this, we define a map $ \eps_i : V \rightarrow \N $ giving the nilpotence degree of $ e_i$ on a vector $ v \in V $; more precisely, $ \eps_i(v) =\max\{n\in\N : e_i^n v\neq0\}$.
	
	A good basis $ B  $ of $ V $ is called \textbf{perfect}, if for each $ i\in I $, and $ b \in B $, either $ e_i b = 0 $ or there exists $ \tilde e_i(b) \in B $ such that
	$$e_i b=\eps_i(b)\;\tilde e_i(b)+ v \quad \text{ for some $ v \in \ker e_i^{\eps_i(b) - 1}$} $$
	
	In other words (up to a predictable scalar)  $ e_i b $ equals $ \tilde e_i(b) $ modulo a vector with lower nilpotence degree.  Note that this definition only requires $ V $ to be a representation of the Borel subalgebra $ \mathfrak b $.
	
	
	\begin{Example}
		To understand these scalars and gain some intuition, it is instructive to consider the case of $ \fg = \ssl_2 $.  In this case, $ P_+ = \N $ and $ V(n) = \C[x,y]_n $, the space of homogenous polynomials of degree $ n$.  The Chevalley generator $ e $ acts by $ y \partial_x $ (on both the left and right) and the unique perfect basis (up to a scalar) is 
		$ \{ x^n,\, x^{n-1} y,\, \dots, y^n \}$.
		
		Note that $ y \partial_x (x^k y^{n-k}) = k x^{k-1} y^{n-k+1} $ and $ \eps(x^k y^{n-k}) = k $.  In this case there is no lower order term. 
	\end{Example}
	
	\begin{Example} \label{eg:sl3adjoint}
		The simplest irreducible representation where lower order terms occur is the adjoint representation of $ \ssl_3$.  In this representation, $ V = \ssl_3 $ with the action given by matrix commutator.  If we assume that $ B $ contains the highest weight vector $ \begin{psmallmatrix} 0 0 1 \\ 0 0 0 \\ 0 0 0 \end{psmallmatrix} $, then it is easy to see that the perfect basis condition forces $B$ to contain 
$$		 \begin{psmallmatrix} 0 & 0 & 1 \\ 0 & 0 & 0 \\ 0 & 0 & 0 \end{psmallmatrix},\, \begin{psmallmatrix} 0 & -1  & 0 \\ 0 & 0 & 0 \\ 0 & 0 &0 \end{psmallmatrix}, \, \begin{psmallmatrix} 0 & 0 & 0 \\ 0 & 0 & -1 \\ 0 & 0 & 0 \end{psmallmatrix},\, \begin{psmallmatrix} 0 & 0 & 0 \\ 1 & 0 & 0 \\ 0 & 0 & 0 \end{psmallmatrix},\, \begin{psmallmatrix} 0 & 0 & 0 \\ 0 & 0 & 0 \\ 0 & 1 & 0 \end{psmallmatrix},\, \begin{psmallmatrix} 0 & 0 & 0 \\ 0 & 0 & 0 \\ -1 & 0 & 0 \end{psmallmatrix}$$
The choice of basis for the diagonal matrices is more interesting.  The requirement that $ B $ be compatible with the kernels of $ e_1, e_2 $ forces $ B $ to contain matrices of the form
$$		 \begin{psmallmatrix} -a & 0 & 0 \\ 0 & -a & 0 \\ 0 & 0 & 2a \end{psmallmatrix},\ \begin{psmallmatrix} 2b & 0 & 0 \\ 0 & -b & 0 \\ 0 & 0 & -b \end{psmallmatrix} $$
for some non-zero $a, b$.  We are forced to take  $  b = 1/3$ and similarly $ a = 1/3$, since
$$e_1 \begin{psmallmatrix} 0 & 0 & 0 \\ 1 & 0 & 0 \\ 0 & 0 & 0 \end{psmallmatrix} = \begin{psmallmatrix} 1 & 0 & 0 \\ 0 & -1 & 0 \\ 0 & 0 & 0 \end{psmallmatrix} = 2 \begin{psmallmatrix} 2/3 & 0 & 0 \\ 0 & -1/3 & 0 \\ 0 & 0 & -1/3 \end{psmallmatrix} + \begin{psmallmatrix} -1/3 & 0 & 0 \\ 0 & -1/3 & 0 \\ 0 & 0 & 2/3 \end{psmallmatrix}
$$
where the second term is of lower nilpotence degree (since it lies in the kernel of $e_1 $).
	\end{Example}
	
	\subsection{Perfect bases and crystals}
	Any perfect basis gives rise to a combinatorial structure called a crystal. 
	Crystals were first introduced by Kashiwara \cite{Kash} as the $q=0 $ limit of a basis for a representation of a quantum group.  However, we prefer to view them as recording the leading order behaviour of $ e_i $ acting on a perfect basis.
	
	 A \textbf{crystal} is a finite set $ B $, along with a map $ \wt : B \rightarrow P $, and for each $ i \in I $, a partially defined map $ \tilde e_i : B \rightarrow B $. If $ B $ is a perfect basis, then it automatically acquires a crystal structure.  The following result of Berenstein--Kazhdan \cite[Theorem~5.37]{BK} shows that this combinatorial structure depends only on the representation. 
	\begin{Theorem} \label{th:uniquerep}
		Let $ V $ be a representation and let $ B, B' $ be two perfect bases.  Then there exists a bijection $ B \cong B' $ which is an isomorphism of crystals.  
	\end{Theorem}
	
	Because of this theorem, we may speak of \textit{the} crystal of a representation.  In particular, the crystal of $ V(\lambda) $ is denoted $ B(\lambda)$.  Many different explicit combinatorial realizations of $ B(\lambda) $ are possible.  In this talk, we will focus on MV polytopes \S \ref{se:MVpol}.

%
%
%

	\subsection{Biperfect bases}
	Rather than looking at each irreducible representation individually, we can study them all at once, using the following trick.  The maximal unipotent subgroup $ N $ has left and right actions of $ N $ and thus the coordinate ring  $ \C[N] $ has left and right actions of $ \mathfrak n $ by differential operators.  For each $ i \in I$, we write $e_i : \C[N] \rightarrow \C[N]$ for the left action and $ e_i^* : \C[N] \rightarrow \C[N] $ for the right action.
	
	For each $ \lambda \in P_+$, choose a highest weight vector $ v_\lambda \in V(\lambda) $ and  let $v_\lambda^*:V(\lambda)\to\mathbb C$ be a dual linear form.  We define an $N$-equivariant map
	$$\Psi_\lambda:V(\lambda)\to\C[N] \qquad  \Psi_\lambda(v)(g)=v_\lambda^*(g v)$$ This linear map is injective and its image is 
\begin{equation} \label{eq:impsi}
	 \im \Psi_\lambda = \bigcap_{i\in I} \ker \, (e_i^*)^{\alpha_i^\vee(\lambda) + 1} \subset \C[N]
	\end{equation}
Thus a basis for $ \C[N] $ compatible with the kernels of all powers of $e_i^* $ gives a basis for each $ V(\lambda) $.  Conversely, a collection of bases for each $ V(\lambda) $ can sometimes glue together to give a basis for $ \C[N]$.


	A basis $ B $ of $ \C[N] $ is called \textbf{biperfect} if it contains $ 1$, and it is perfect with respect to both the left and right actions of $ \mathfrak n $.  Thus, $ B $ will have two families of crystal operators, written $ \tilde e_i, \tilde e^*_i $ and two families of maps $ \eps_i, \eps_i^* : B \rightarrow \N$ (but only one weight map).  
	
		From Proposition \ref{pr:Zhelo} and (\ref{eq:impsi}), we immediately deduce the following corollary, which can be regarded as a generalization (from the canonical basis to arbitrary biperfect bases) of \cite[Corollary 3.4]{BZ}.
	\begin{Corollary} \label{co:biperfectrestrict}
		Let $ B $ be a biperfect basis of $ \C[N] $.  
		\begin{enumerate}
			\item For any $ \lambda \in P_+$, the set
			$ \{ b \in B : \eps_i^*(b) \le \alpha_i^\vee(\lambda) \} 
			$ restricts via $\Psi_\lambda $ to a perfect basis for $ V(\lambda) $.
		\item
		For any $ \lambda, \mu, \nu \in P_+ $, the set
		$$
		\{b\in B : 
		\wt(b)=\nu-\mu - \lambda\ \text{ and }\ \forall i\in I,\ \varepsilon_i(b)
		\leq \alpha_i^\vee(\mu), \ \varepsilon_i^*(b) \leq \alpha_i^\vee(\lambda) \} 
		$$
		restricts to a basis for $C^\nu_{\lambda \mu} $.
		\end{enumerate}
	\end{Corollary}
	Thus we can solve Problem \ref{prob:mult} by understanding well the bicrystal structure on $ B$.

\subsection{The bicrystal $B(\infty) $}	 \label{se:Binf}
	The Berenstein--Kazhdan result (Theorem \ref{th:uniquerep}) generalizes to biperfect bases.
	\begin{Theorem} \cite[Theorem 2.4]{mvbasis}
		\label{th:unique}
		Let $B$ and $B'$ be two biperfect bases of $\C[N]$. Then there is a unique
		bijection $B\cong B'$ that respects the bicrystal structure.
	\end{Theorem}

The abstract combinatorial crystal underlying any biperfect basis is denoted $ B(\infty) $.  On $B(\infty) $ we have the Kashiwara involution $ * : B(\infty) \rightarrow B(\infty) $ exchanging $ \tilde e_i $ and $ \tilde e_i^* $.

	\begin{Remark}
	The algebra $ \C[N] $ has an involutive automorphism $*$ (coming from the inverse map on $ N $) which exchanges the left and right actions of $ \fn $.  The Kashiwara involution $ *$ on $ B(\infty)$ is the combinatorial manifestation of the involution $* $ on $ \C[N]$.
	
	 A $ *$-\textbf{invariant perfect basis} is a perfect basis which is invariant under $ *$.  Every example of a biperfect basis that we know is $*$-invariant.  This was established by Lusztig \cite[Prop. 3.3]{Lusztig} for the dual canonical basis (when $ \fg $ is simply-laced); in fact this was the first appearance of the Kashiwara involution in the literature. \end{Remark}

\begin{Example}When $ G = SL_3$, $B(\infty) $ can be drawn in the following way. The action of $ e_1 $ is given by right-pointing diagonal arrows and the action of $ e_2 $ is given by the left-pointing ones.  Each horizontal group of dots have the same weight and the Kashiwara involution flips each such group. We would like to thank Mark Haiman for showing us this drawing many years ago.
	
	
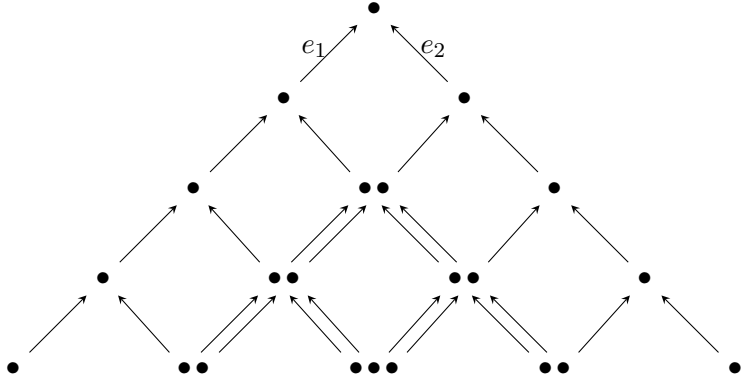
\begin{center}
	\begin{tikzpicture}[scale=1.2]
		\node (001) at (0,0){$\bullet$};
		\node (011) at (-1,-1){$\bullet$};
		\node (101) at (1,-1){$\bullet$};
		
		\node (021) at (-2,-2){$\bullet$};
		\node (111) at (-0.1,-2){$\bullet$};
		\node (112) at (0.1,-2){$\bullet$};
		\node (201) at (2,-2){$\bullet$};
		
		\node (031) at (-3,-3){$\bullet$};
		\node (121) at (-1.1,-3){$\bullet$};
		\node (122) at (-0.9,-3){$\bullet$};
		\node (211) at (0.9,-3){$\bullet$};
		\node (212) at (1.1,-3){$\bullet$};	
		\node (301) at (3,-3){$\bullet$};
		
		\node (041) at (-4,-4){$\bullet$};
		\node (131) at (-2.1,-4){$\bullet$};
		\node (132) at (-1.9,-4){$\bullet$};
		\node (221) at (-0.2,-4){$\bullet$};
		\node (222) at (0,-4){$\bullet$};	
		\node (223) at (0.2,-4){$\bullet$};
		\node (311) at (1.9,-4){$\bullet$};	
		\node (312) at (2.1,-4){$\bullet$};	
		\node (401) at (4,-4){$\bullet$};
	
		\draw[-stealth] (011) -- node[near start,above] {$e_1$} (001);		
		\draw[-stealth] (101) -- node[near start,above] {$e_2$}  (001);
		
		\draw[-stealth] (021) to (011);
		\draw[-stealth] (111) to (011);
		\draw[-stealth] (112) to (101);
		\draw[-stealth] (201) to (101);
	
		\draw[-stealth] (031) to (021);
		\draw[-stealth] (121) to (021);
		\draw[-stealth] (121) to (111);
		\draw[-stealth] (122) to (112);
		\draw[-stealth] (211) to (111);
		\draw[-stealth] (212) to (112);
		\draw[-stealth] (212) to (201);
		\draw[-stealth] (301) to (201);
		
		\draw[-stealth] (041) to (031);
		\draw[-stealth] (131) to (031);
\draw[-stealth] (131) to (121);
\draw[-stealth] (132) to (122);
\draw[-stealth] (221) to (121);
\draw[-stealth] (222) to (122);
\draw[-stealth] (222) to (211);
\draw[-stealth] (223) to (212);
\draw[-stealth] (311) to (211);
\draw[-stealth] (312) to (212);
\draw[-stealth] (312) to (301);
\draw[-stealth] (401) to (301);
		
\end{tikzpicture}
\end{center}
\end{Example}

\subsection{Biperfect bases in small rank}
	For small rank groups, it is easy to show the existence and uniqueness of biperfect bases of $\C[N] $ by elementary means.
	
	\begin{Theorem} \label{th:uniquebasis}
	For $ G=\SL_2, \SL_3, \SL_4$, $\C[N] $ has a unique biperfect basis.
	\end{Theorem}
	
	\begin{Example}
		Suppose $ G = \SL_2$, then $ \C[N] = \C[x] $ and $ \Psi_n : \C[x,y]_n \rightarrow \C[x] $ is the map sending $ y $ to $ 1$.  The left and right actions of $ e \in \mathfrak n $ on $\C[x] $ agree and are given by $ e = \partial_x $. The unique biperfect basis of $ \C[x] $ is $ \{ 1, x, x^2, \dots \} $.
	\end{Example}
	\begin{Example}
		Suppose $G=\SL_3$, with the standard choice for $B$, $T$ and
		$N$. Then $\C[N]=\mathbb C[x,y,z]$ where $x$, $y$ and $z$ are the three
		matrix entries of an upper unitriangular matrix
		$$\begin{psmallmatrix}1&x&z\\0&1&y\\0&0&1\end{psmallmatrix}\in N.$$
		The unique biperfect basis of $\C[N]$ is
		$$B=\{x^az^b(xy-z)^c : (a,b,c)\in \N^3\}\cup
		\{y^az^b(xy-z)^c : (a,b,c)\in \N^3\}.$$
	\end{Example}

\subsection{Three different biperfect bases}
	
	For general $ G$, biperfect bases are not unique, nor is it very easy to show their existence.  
	
	The first example of a biperfect basis was Lusztig's \textbf{dual canonical basis} which is also known as Kashiwara's upper global basis \cite{Lusztig90,Kashiwara93a}.  This is actually a basis for the corresponding quantum deformation of $ \C[N] $, but it can be specialized at $ q = 1$ to give a biperfect basis.
	
	Another example, when $ \fg $ is simply-laced, is Lusztig's \textbf{dual semicanonical basis} \cite{Lusztig00}, which is constructed by means of the representation theory of the preprojective algebra.
	
	A third example is the \textbf{Mirkovi\'c--Vilonen basis} \cite{MV} coming from the geometry of the affine Grassmannian.
	
	This trichotomy of bases will be the focus of this talk.  Each of these bases comes from a complicated algebraic or geometric source.  Following Arun Ram, we can imagine three high mountains whose springs give these three bases.
	
	These bases are all different, combining \cite[Thm 1.7]{mvbasis}, \cite[Prop 2.7]{BGL}, \cite[(3)]{GLSsemi}.
\begin{Theorem}  \label{th:counterexample}
	For $ G = \SL_6 $ in the weight space $2\alpha_1 + 4\alpha_2 + 4\alpha_3 + 4\alpha_4 + 2\alpha_5$, and for $ G = \SO_8$ in the weight space $ 2\alpha_1 + 4 \alpha_2 + 2\alpha_3 + 2 \alpha_4$, there is a point of $B(\infty) $ whose corresponding dual canonical, dual semicanonical, and MV bases are all different.
	
	Moreover, in both these examples, we have the following specific situation.
\begin{equation}	\label{eq:bcd}
	d = b + v \qquad c = b + 2v
\end{equation}
	where $b, c, d$ denote the MV, dual semicanonical, and dual canonical basis vectors, all of whom define the same point in $ B(\infty) $, and $ v $ denotes a vector common to all three bases.
\end{Theorem}


\begin{Question} \label{qu:biperfectset}
	What can we say about the set of all biperfect bases of $ \C[N] $ for a fixed $ G$?
\end{Question}
	
\section{Mirkovi\'c--Vilonen basis}

\subsection{MV cycles}
Mirkovi\'c--Vilonen \cite{MV} used the geometric Satake correspondence to define the MV basis for irreducible representations of $ G $.  This basis is indexed by certain subvarieties in the \textbf{affine Grassmannian}, known as MV cycles.

Let $ G^\vee $ be the Langlands dual group and let
$ \Gr = G^\vee\! (\!( t )\!) /G^\vee \llbracket t \rrbracket $ denote the affine Grassmannian of
this group.  By definition, the coweight lattice of $ G^\vee $ coincides with the weight lattice $P $ of $ G $.  For each coweight $ \mu \in P$, we get a point of $G^\vee\! (\!(t)\!) $ and hence a point $ L_\mu $ in $ \Gr $.  
Let $ S_{\pm}^\mu := N^\vee_{\pm}\! (\!(t)\!) L_\mu $ denote
semi-infinite orbits in $ \Gr$, where $ N^\vee_{\pm} $ denote opposite unipotent subgroups in $ G^\vee $.

For $\lambda \in P_+ $, let $ \Gr^\lambda := \overline{G^\vee\llbracket t\rrbracket L_\lambda} $ be a \textbf{spherical Schubert variety}.  This is a finite-dimensional singular projective variety whose geometry is closely related to the irreducible representation $ V(\lambda) $.  Let $ P(\Gr)$ denote the category of perverse sheaves on $ \Gr $ which are constructible with respect to the stratification by $ G^\vee\llbracket t\rrbracket$ orbits.  This is a semisimple category whose simple objects are the intersection cohomology sheaves $ IC_\lambda $ of the spherical Schubert varieties.  There is a monoidal structure on $ P(\Gr) $ by convolution.

The following geometric Satake correspondence was established by Mirkovi\'c--Vilonen \cite{MV}, following earlier work by Lusztig \cite{Lusztig83} and Ginzburg \cite{Ginzburg}.
\begin{Theorem} \label{th:Satake}
	\begin{enumerate}
\item  There is an equivalence of monoidal categories, $ P(\Gr) \cong \Rep G $.
  \item Under this equivalence, for each $ \lambda \in P_+$, $ IC_\lambda $ is sent to $ V(\lambda) $. 
  \item Under this equivalance, for each $ \mu \in P$, the hyperbolic stalk functor $ H^\bullet_{S_-^\mu}(-) $ matches the functor of taking the $ \mu$-weight space. 
  \end{enumerate}
 \end{Theorem}
 Combining these statements we conclude
$H_{\ttop}(\Gr^\lambda \cap S_-^\mu) \cong V(\lambda)_\mu $.
The irreducible components of $ \Gr^\lambda \cap S_-^\mu $ are called \textbf{MV cycles}.  Via this theorem, they provide a basis for each $ V(\lambda)_\mu $.

\subsection{Stable MV cycles}
For bases of $ \C[N] $, we will be concerned with the
intersection of opposite semi-infinite orbits.  For any
$ \nu \in Q_+ $, the positive root cone, the irreducible
components of $ \overline{S_+^\nu \cap S_-^0} $ are called
\textbf{stable MV cycles}.

Given an MV cycle  $ Z \subset \Gr^\lambda \cap S_-^\mu $, we can translate by $ t^{-\mu} $ to produce a stable MV cycle $ \overline{t^{-\mu} Z } $.  This process is the geometric analog of the map $ \Psi_\lambda : V(\lambda) \rightarrow \C[N] $.

In \cite{mvbasis}, we combined work of Ginzburg \cite{Ginzburg} and Mirkovi\'c--Vilonen \cite{MV} to prove the following result, which had been conjectured by Anderson \cite{Anderson03}.  

\begin{Theorem}
\begin{enumerate}
	\item The MV bases for each $ V(\lambda) $ can be collected together to form a biperfect basis for $\C[N]$, which is indexed by stable MV cycles.
	\item For each $ i $, the action of $ e_i $ on an MV basis vector $ b_Z $ is given by the intersection of the stable MV cycle $ Z $ with a hyperplane section.
	\item Given two MV cycles $ Z_1, Z_2 $, the product $ b_{Z_1} b_{Z_2} $ in $ \C[N] $ is given by the Beilinson--Drinfeld degeneration of $ Z_1 \times Z_2 $.
\end{enumerate}
In particular, the structure constants for the action of $ e_i $ and for the multiplication are non-negative integers.
\end{Theorem}

\subsection{MV polytopes} \label{se:MVpol}
For each stable MV cycle $ Z$, we define $ \Pol(Z) $ to be its moment map image (for a real Hamiltonian torus action).  Equivalently, we have
$$ \Pol(Z) = \Conv \bigl(\mu : L_\mu \in Z \bigr) $$
The polytopes produced this way are called MV polytopes.  In \cite{mvcycle}, we proved the following result.

\begin{Theorem} \label{th:2faces}
	The map $Z \mapsto \Pol(Z) $ gives a bijection between the stable MV cycles and the MV polytopes.  The MV polytopes are precisely those lattice polytopes whose dual fan is a coarsening of the Weyl fan and whose 2-faces are MV polygons for the appropriate rank 2 groups (which can be described explicitly).
\end{Theorem}

This theorem was reinterpreted by Goncharov--Shen \cite{GS} as the following statement.

\begin{Corollary} \label{co:GS}
	The MV polytopes are in natural bijection with $(G^\vee/B^\vee)(\Z_{\mathrm{trop}})_\ge $, the non-negative tropical points of the flag variety.
\end{Corollary}

	Following the historical order, we have described MV polytopes as the moment map images of MV cycles.  However, we emphasize that Theorem \ref{th:2faces} shows that they are purely combinatorial objects.  We will see in the next two sections that these same polytopes are naturally obtained from general preprojective algebras modules and simple KLR modules.  They are the common shadows from all three mountains.
	
	In \cite{mvcrystal}, we gave an explicit description of the crystal structure on the set of MV polytopes.  This provides a convenient combinatorial framework for describing the crystal $ B(\infty) $ and is easily connected to many other combinatorial models.  In particular, for each reduced word $ s_{i_1} \cdots s_{i_m} = w_0$ for the longest element of the Weyl group, Lusztig \cite{Lusztig} constructed a bijection $ B(\infty) \rightarrow \N^m $ using the relation between PBW monomials and the canonical basis.  In \cite{mvcrystal}, we showed that the Lusztig datum of  $ b \in B(\infty) $ is the list of lengths along a path following the edges of the MV polytope $ \Pol(b) $, in root directions determined by the reduced word.

\begin{Example}
	Take $ G = \SL_3$.  In this case, an MV polytope is a hexagon with all $120^\circ$ angles, whose ``width'' $A$ is equal to the maximum of its two ``heights'' $B, C$.
	
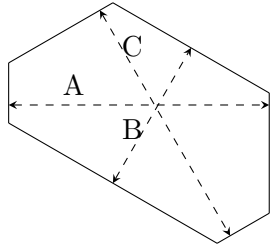
	\begin{center}
	\begin{tikzpicture}[scale=0.8]
		\coordinate (e) at (0,0){};
		\coordinate (s1) at ({-sqrt(3)},-1){};
		\coordinate (s2) at ({3/2*sqrt(3)},-1.5){};
		\coordinate (s1s2) at ({-sqrt(3)},-2){};
		\coordinate (s2s1) at ({3/2*sqrt(3)},-3.5){};
		\coordinate (w0) at ({sqrt(3)},-4){};
		\draw (e) to (s1);		
		\draw (e) to (s2);		
		\draw (s1) to (s1s2);		
		\draw (s2) to (s2s1);		
		\draw (w0) to (s1s2);		
		\draw (w0) to (s2s1);
		\draw [stealth-stealth,dashed]({-sqrt(3)},-1.7) --  node[above,near start] {A}
  ({3/2*sqrt(3)},-1.7);
  		\draw [stealth-stealth,dashed](0,-3) --  node[above,near start] {B}
  ({3/4*sqrt(3)},-0.75);	
    		\draw [stealth-stealth,dashed]({-1/8*sqrt(3)},{-1/8}) --  node[above,near start] {C}
  ({9/8*sqrt(3)},{-31/8});	
	\end{tikzpicture}
	\end{center}
	For this polytope, the two Lusztig data are $(3,2,1) $ and $ (2,1,4) $.
\end{Example}

\section{Dual semicanonical basis}
\subsection{Preprojective algebra}
Assume for this section that $ \fg $ is simply-laced.  
 Let $ H $ denote the set of oriented edges of the Dynkin diagram of $ \fg $.  If $ h = (i,j) $, write $ \bar{h} = (j,i) $.  Fix a map $ \tau : H \rightarrow \{1,-1\} $ such that for each $ h $, $ \tau(h) + \tau(\bar{h}) = 0 $ (such a $\tau $ corresponds to an orientation of each edge of the Dynkin diagram).

The \textbf{preprojective algebra} $ \Lambda $ is the quotient of the
path algebra of $ (I,H) $ by the relation $\sum_{h \in H} \tau(h) h \bar{h} = 0 $. So a $ \Lambda$-module $ M $ consists of vector spaces $ M_i $, for $ i \in I $, and linear maps $ M_h : M_i \rightarrow M_j $ for each $ h = (i,j) \in H $, such that
\begin{equation} \label{eq:preproj}
	\sum_{h \in H} \tau(h) M_h M_{\bar{h}} = 0 
\end{equation}
Given a $ \Lambda$-module $ M$, we define its dimension vector by
$$
\dimvec M = \sum_{i \in I} (\dim M_i)\, \alpha_i
$$

We write $ S_i $ for the simple module at vertex $ i $, the unique module with $ \dimvec S_i = \alpha_i $.

For each $ \nu = \sum_{i \in I} \nu_i \alpha_i  \in Q_+ $, we consider the affine variety of $ \Lambda$-module structures on $ \bigoplus_{i \in I} \C^{\nu_i} $.  More precisely, we define
$$
\Lambda(\nu) \subset \bigoplus_{(i,j) \in H} \Hom(\C^{\nu_i}, \C^{\nu_j})
$$
to be the subvariety defined by equation (\ref{eq:preproj}).

\subsection{The dual semicanonical basis}\label{ssec:dualsemi}
Let $ M$ be a $ \Lambda $-module.
Following Lusztig \cite{Lusztig00} and Geiss--Leclerc--Schr\"oer \cite[\S 5]{GLSsemi}, we define an element
$ \xi_M \in \C[N]$ as follows.  First, for each $ \vi \in I^p $, we define the projective variety of composition series of type $ \vi$,
$$
F_\vi(M) = \bigl\{ 0=M^0 \subset M^1 \subset \dots \subset M^p = M : M^k/M^{k-1} \cong S_{i_k} \text{ for all } k \bigr\}
$$
and then we define $\xi_M \in \C[N]$ by requiring that
$$ \langle e_{i_1} \cdots e_{i_p} , \xi_M \rangle = \chi(F_\vi(M))$$
for any $\vi \in I^p$, where $ \chi $ denotes topological Euler characteristic, and where $ \langle \, , \rangle $ denotes the pairing between $ U \mathfrak n $ and $ \C[N] $.

This map $ M \mapsto \xi_M $ is constructible and so for any irreducible component $ Y \subset \Lambda(\nu) $, we can define $ c_Y \in \C[N]_{\nu} $ by setting $ c_Y = \xi_M $, for $ M $ a general point in $ Y $.

The following result is due to Lusztig \cite{Lusztig00}.
\begin{Theorem} \label{th:LusztigPerfect}
	\begin{enumerate}
		\item
		For each $ \nu \in Q_+$, $ \{ c_Y \mid Y \in \Irr \Lambda(\nu) \} $ is a basis for~$ \C[N]_{\nu}$.
		\item
		Together they form a biperfect basis of $ \C[N] $, called the dual semicanonical basis.
	\end{enumerate}
\end{Theorem}

\subsection{Polytopes from preprojective algebra modules}

In the resulting bicrystal structure on the set $ \sqcup_\nu \Irr \Lambda(\nu) $, we have
$$
\eps_i(Y) = \dim \Hom_\Lambda(M, S_i) \quad \eps_i^*(Y) = \dim \Hom_\Lambda(S_i, M)
$$
where $ M $ is a general point of $ Y$.  

\begin{Remark}
Fix $ \lambda \in P_+$.  From Corollary \ref{co:biperfectrestrict}(1), those components $ Y $ which satisfy $ \eps_i^*(Y) \le \alpha_i^*(\lambda) $ index a basis for $ V(\lambda) $.  These same irreducible components form the core of the corresponding Nakajima quiver variety via the correspondence explained in \cite[Section 4.6]{Saito}.
\end{Remark}

The bicrystal $\sqcup \Irr \Lambda(\nu) $ is isomorphic to $ B(\infty) $ by Theorem \ref{th:unique}.  Thus, a MV polytope is canonically associated to component $ Y$.  We can describe this polytope using the module structure on a general point of $ Y$.

\begin{Theorem} \cite[\S 1.3]{BKT} \label{thm:BKT}
	Let $ Y $ be a component of $ \Lambda(\nu) $ and $M $ be a general point of $ Y$.  The  MV polytope of the basis vector $c_Y$ is given by the \textbf{Harder--Narasimhan} polytope of $ M$.
$$
\Pol(M) := \Conv \bigl(\dimvec N : N \subseteq M \text{ is a submodule} \bigr)
$$	
\end{Theorem}

\begin{Example}
	Take $ \fg  = \ssl_3$, and take $\nu = \alpha_1 + \alpha_2 $.
	
	Then $ \Lambda(\nu) = \{(a,b) \in \C^2 : ab = 0 \} $ where $(a,b)$ corresponds to the $\Lambda$-module
	$$ \C \stackrel[a]{b}{\leftrightarrows} \C $$
	$\Lambda(\nu) $ has two components in this case.  If $ Y $ is the component given by $ b= 0 $, then for general $M \in Y $, we have submodules in $ M $ of dimension $ 0, \alpha_1, \alpha_1 + \alpha_2$ and so $ \Pol(M) $ is the triangle with these vertices.
\end{Example}
	
\section{Dual canonical bases}
\subsection{KLR algebras}
Let $ \fg $ an arbitrary semisimple Lie algebra.  For each $ \nu \in Q_+ $, Khovanov--Lauda--Rouquier defined an algebra $ R_\nu $.  It can either be defined by generators and relations using a modification of the presentation of a Hecke algebra \cite{R}, or as an algebra of decorated string diagrams as in \cite{KL}.  This algebra can also be realized as an Ext algebra of certain perverse sheaves constructed by Lusztig (see \cite{VV}). 

In this section, we work over $ \C$, though it is possible to work over fields of positive characteristic; this produces different bases, called the dual $p$-canonical bases.

The algebra $ R_\nu $ contains indempotents $ e_\vi $ indexed by sequences $ (i_1, \dots, i_p) \in I^p $ such that $ \alpha_{i_1} + \cdots + \alpha_{i_p} = \nu$.   We write $ K(R_\nu)_\C$ for the complexified Grothendieck group of finite-dimensional $ R_\nu $ modules. The following result is due independently to Rouquier and Khovanov--Lauda.

\begin{Theorem} \label{th:KLR}
	For each $ \nu \in Q_+ $, there is an isomorphism $ K(R_\nu)_\C \cong \C[N]_{\nu} $, written $ [L] \mapsto d_L $ such that
	$ \langle e_{\vi}, d_L \rangle = \dim e_\vi L $,
	for any module $ L $ and $ \vi $ as above.
\end{Theorem}

The vector space $ K(R_\nu)_\C$ has a basis given by the simple finite-dimensional $R_\nu$-modules.  By \cite{VV}, the resulting basis for $ \C[N] $ coincides with Lusztig's dual canonical basis.

\subsection{Polytopes from KLR modules}
  The bicrystal structure on the set of simple KLR modules was carefully studied by Lauda--Vazirani \cite{LV}.  By Theorem \ref{th:unique}, this bicrystal is isomorphic to $ B(\infty) $, and thus an MV polytope is canonically associated to each simple KLR module.

On the other hand, these algebras come with non-unital morphisms $ R_\mu \otimes R_{\nu - \mu} \rightarrow R_\nu $. We write $ e_{\mu, \nu - \mu} $ for the image of the identity under this map.  Tingley--Webster \cite{TW} used these morphisms to define a polytope associated to each KLR module.

\begin{Theorem}
	The MV polytope of a simple $ R_\nu $-module $L $ is the \textbf{character polytope}
	$$ \Pol(L) := \Conv \bigl(\mu : e_{\mu, \nu - \mu}L \ne 0 \bigr) $$
\end{Theorem}

\subsection{Generalizations to affine and Kac--Moody cases}
Unlike the MV basis, the dual semicanonical basis and dual canonical basis admit straightforward generalizations to the setting where $ \fg $ is a symmetric (resp. symmetrizable) Kac--Moody Lie algebra.

The polytopes $ \Pol(M) $ and $ \Pol(L) $ associated to a general $ \Lambda$-module or a simple KLR module admit obvious generalizations in this setting.  However, due to higher root multiplicities, there can't be a Lusztig data bijection $ B(\infty) \rightarrow \N^{\Delta_+} $ (where $ \Delta_+ $ is the set of positive roots).  Thus, we must enhance the polytope with some extra information.  This was carried out in \cite{BKT} (using $\Lambda$-modules) and in \cite{TW} (using KLR modules).

In the affine case, this extra information consists of partitions associated to vertical edges of the polytope (vertical edges are those pointing in the imaginary root direction).  Moreover, these decorated polytopes are characterized by their 2-faces (as in the finite case, Theorem \ref{th:2faces}) and the new relevant polygons were described combinatorially in \cite{BDKT}.

\begin{Question}
	\begin{enumerate}
\item	Is it possible to give a ``tropical'' description of affine MV polytopes, similar to Corollary \ref{co:GS}?
\item
	Outside of the affine type, it is possible to give a combinatorial description to the extra information carried on the MV polytope?

	\item Though more complicated, the theory of MV cycles exists for affine Kac--Moody Lie algebras.  How can we relate these MV cycles to the affine MV polytopes?  In particular, what information about the cycles is encoded in the partitions along vertical edges?
	\end{enumerate}
\end{Question}

\section{Comparing biperfect bases}
\subsection{Change of basis matrix}
Let $ B, B' $ be two biperfect bases for $ \C[N] $.  By Theorem \ref{th:unique}, we obtain bijections $ B \rightarrow B(\infty) \leftarrow B' $ and thus a bijection between $ B $ and $ B'$.  Thus, it makes sense to speak of the change of basis matrix between $ B $ and $ B'$.  In \cite{Baumann}, Baumann proved that this matrix is upper unitriangular with respect to a partial order on $ B(\infty) $, defined combinatorially using the crystal structure.  Many elements of $ B(\infty)$ are incomparable using this order.  Thus, many off-diagonal elements of the change of basis matrix must vanish.  In low rank, the order becomes trivial and gives the proof of Theorem \ref{th:uniquebasis}.

\subsection{Measures} \label{se:measures}
We think of an MV polytope as the shadow of a biperfect basis vector.  Unfortunately, this shadow is not precise enough to distinguish between different biperfect basis vectors which represent the same element of $B(\infty) $.  For this reason, we now introduce a measure supported on the MV polytope.

Consider the vector space $\Dist(\tR) $ of $\C$-valued compactly supported
distributions on $ \tR $.  It forms an algebra under convolution, using the addition map $\tR \times \tR \xrightarrow{+} \tR$.

Let $ \Delta^p := \{(c_0, \dots, c_p) \in \mathbb R^{p+1} : \text{ each }c_i\geq0,\ c_0 + \cdots + c_p = 1 \} $ be the standard $p$-simplex.
For $ \vi \in I^p $, we define the linear map
$ \pi_\vi : \R^{p+1} \rightarrow \tR $ by
$$ \pi_\vi(c_0,\dots,c_p) = \sum_{k = 0}^p c_k\,(\alpha_{i_1} + \cdots + \alpha_{i_k})$$

We define the measure $ D_\vi $ on $\tR$ by $ D_\vi := (\pi_\vi)_*(\delta_{\Delta^p}) $, the push-forward of Lebesgue measure on the $p$-simplex.  	The measures $ D_\vi $ satisfy the shuffle identity.

\begin{Lemma} \cite[Lemma 8.5]{mvbasis} \label{le:shuffle}
For $ \vj \in I^p, \vk \in I^q$, 
	$$
	D_\vj * D_\vk = \sum_{ \vi \in \vj \shin \vk} D_\vi.
	$$
	where $ \vj \shin \vk $ is the set of all sequences obtained by shuffling $ \vj $ and $ \vk $.
\end{Lemma}

Elementary considerations involving the coproduct structure on $ U \fn $ show that the shuffle identity implies that there is an algebra morphism $ D : \C[N] \rightarrow \Dist(\tR) $ defined by
$$
D(f) = \sum_{\vi} \langle e_{\vi}, f \rangle D_\vi
$$

\subsection{Fourier transform} \label{se:FTmeasures}
For each weight $\beta\in P$, we define $e^\beta$ to be the function $x\mapsto e^{\langle \beta, x \rangle} $ on $ \ft_\C $.  Given a measure $ D(f) $ as above, we can consider its Fourier transform $ FT(D(f)) $ which lies in the space of meromorphic functions on $ \ft_\C $, spanned by these exponentials over the field $\C(\ft)$ of rational functions. In this way, we obtain the following result.

\begin{Proposition} \label{pr:FTD}
	The composition $ FT \circ D $ defines an algebra morphism
	$$
	\C[N] \rightarrow \C(\ft) \otimes \C[T]
	$$
\end{Proposition}
This algebra morphism defines a rational map $ \ft \times T \rightarrow N $. This map is actually regular on $ \ftreg \times T  $, where $ \ftreg $ is the complement of the root hyperplanes in $ \ft $. 
\begin{Theorem}\label{th:FTDpsi} \cite[Theorem 8.11]{mvbasis}
	\begin{enumerate}
		\item For all $ x \in \ftreg$, there exists a unique $ n_x \in N $ such that $ \Ad_{n_x}(x)=x+e$.
		\item The rational map from Proposition \ref{pr:FTD} is given by $ (x,t) \mapsto t^{-1} n_x t n_x^{-1}$.
	\end{enumerate}
\end{Theorem}

We now study a simpler invariant $ \barD $. For a sequence $\vi=(i_1,\ldots,i_p)$, we define
$$\barD_\vi=\prod_{k=1}^{p}\frac{1}{\alpha_{i_k} + \cdots + \alpha_{i_p}} \in \C(\ft)$$
	
\begin{Proposition} \cite[Prop~8.4 and Lemma~8.7]{mvbasis}
	\begin{enumerate}
		\item
These rational functions satisfy the shuffle identity from Lemma \ref{le:shuffle} and thus define an algebra morphism $ \barD : \C[N] \rightarrow \C(\ft) $ by
$
\barD(f) = \sum_{\vi} \langle e_{\vi}, f \rangle \barD_{\vi}
$
\item For any $ f \in \C[N]$, $ \barD(f) $ is the coefficient of $ e^0 $ in $FT(D(f)) $.
\item The algebra morphism from (1) comes from the morphism of varieties $ \ftreg \rightarrow N $ given by $ x \mapsto n_x $.
\end{enumerate}
\end{Proposition}

\subsection{Duistermaat--Heckman measure}
In symplectic geometry, the Duistermaat--Heckman (DH) measure of a symplectic manifold with a Hamiltonian torus action is defined to be the push-forward of the Liouville measure under the moment map.  Brion--Procesi \cite{BrionProcesi} reformulated this notion in algebraic geometry by considering the asymptotics of sections of equivariant line bundles.

Fix a $W$-invariant bilinear form on $ \ft $ (normalized so that short roots have length 1).  This leads to a central extension of $ G^\vee\!(\!(t)\!) $ and thus an equivariant line bundle $ \mathcal O(1) $ on $ \Gr $.

Let $ Z \subset \Gr $ be an MV cycle. The torus $ T^\vee $ acts on the space of sections $ \Gamma(Z, \mathcal O(n))$.  We consider $ [\Gamma(Z, \mathcal O(n))] $ as a distribution on $ \tR $ by 
$$
 [\Gamma(Z, \mathcal O(n))]  = \sum_{\mu \in P} \dim \Gamma(Z, \mathcal O(n))_\mu \, \delta_\mu
$$
(Implicitly, we use the bilinear form to identify $ \ft $ and $ \ft^*$.)

Let $ \tau_n : \tR \rightarrow \tR $ be the automorphism given by scaling by $ \frac{1}{n}$. The \textbf{Duistermaat--Heckman measure} of $Z $
	is defined to be the limit $$ DH(Z) := \displaystyle{\lim_{n \to \infty}} \frac{1}{n^{\dim Z}} (\tau_n)_*[\Gamma(Z, \mathcal O(n))] $$
	within the space of distributions on $ \tR $. Note that each $ (\tau_n)_*[\Gamma(Z, \mathcal O(n))] $ is supported on $ \Pol(Z) $, and hence so is $ DH(Z)$.

Via the Fourier transform, $ DH(Z) $ is  closely related to the class of $ Z $ in the equivariant homology of the affine Grassmannian.  The following ideas are not specific to the affine Grassmannian: they apply to any (ind-) projective variety equipped with a torus action having isolated fixed points.  First, we have the localization theorem in equivariant homology (see for example \cite{Br}).

\begin{Theorem} \label{th:HTGrT}
	The inclusion $\Gr^{T^\vee} \rightarrow \Gr $ induces an isomorphism
	$$
	H^{T^\vee}_\bullet(\Gr^{T^\vee}) \otimes_{\C[\ft]} \C(\ft) \xrightarrow{\sim} H^{T^\vee}_\bullet(\Gr) \otimes_{\C[\ft]} \C(\ft) .$$
\end{Theorem}
Because of this theorem and using $ \Gr^{T^\vee} = \{ L_\mu : \mu \in P \} $, we can write
$$
[Z] = \sum_{\mu \in P} m_\mu(Z)[\{L_\mu\}]
$$
for unique $ m_\mu(Z) \in \C(\ft) $.  Following Brion \cite{Br}, we
call $ m_\mu(Z)$ the \textbf{equivariant multiplicity} of $ Z $
at $L_\mu $.  In \cite{mvbasis}, we proved the following result (again in the general context of projective varieties with torus actions), following ideas of Brion.

\begin{Theorem} \label{th:FTDH} For any MV cycle $ Z$,
	$$FT(DH(Z)) = \sum_{\mu \in P} m_\mu(Z) e^{\mu}
	$$
\end{Theorem}

\subsection{DH measures and measures from $ \C[N] $}

The $T^\vee$-equivariant homology of $ \Gr $ was computed by Yun--Zhu \cite{YunZhu}.  They defined a commutative convolution algebra structure on $ H_\bullet^{{T^\vee}}(\Gr) $ and described this algebra using the geometric Satake correspondence.

	Let $ e = \sum e_i $ be a regular nilpotent element.  We define the \textbf{universal centralizer space} to be
	$$
	C := \bigl\{ (x, b) \in \ft \times B : \Ad_b(x+e) = x +e \bigr\}
	$$
	\begin{Remark}
	For any $ x \in \ft$, $ x+e $ is regular and has centralizer contained in $ B$.  Thus our space $ C$ is the base change over $ \ft \rightarrow \ft/W $ of the usual universal centralizer (often denoted $ J $), as defined in, for example \cite[\S 2.2]{BFM}.
	\end{Remark}
	
	From the definition, we have a map $ C \rightarrow \ft \times T $ given by $ (x,tn) \mapsto (x,t) $. The dual algebra map fits into the following diagram.
\begin{Theorem} \cite[Prop 3.3 and 5.7]{YunZhu} \label{th:HTGr}
	There is an isomorphism of algebras $\theta : \C[C] \rightarrow H_\bullet^{{T^\vee}}(\Gr) $ making the following diagram commute
	$$
	\xymatrix{  
		\C[\ft] \otimes \C[T] \ar^\sim[r] \ar[d]  & H_\bullet^{{T^\vee}}(\Gr^{{T^\vee}}) \ar[d]^{\text{Theorem \ref{th:HTGrT}}} \\
		\C[C] \ar^\theta[r] & H_\bullet^{{T^\vee}}(\Gr)
	}
	$$
\end{Theorem}

Recall the map $ D : \C[N] \rightarrow \Dist(\tR) $ defined in \S \ref{se:measures} and the map $ \barD : \C[N] \rightarrow \C(\ft) $ defined in \S \ref{se:FTmeasures}.  Combining Theorems \ref{th:FTDpsi}, \ref{th:FTDH} and \ref{th:HTGr} we proved the following in \cite{mvbasis}.
\begin{Corollary} \label{co:DDH}
	For any stable MV cycle $Z$, $$ D(b_Z) = DH(Z) \qquad \barD(b_Z) =m_0(Z)$$
\end{Corollary}
	
This Corollary is very useful, since the equivariant multiplicity $ m_0(Z) $ is easily computed using computer algebra programs.  In the appendix of \cite{mvbasis} (written with C. Morton-Ferguson and A. Dranowski), we used this approach to establish part of Theorem \ref{th:counterexample}.

%
%

\subsection{Measures from preprojective algebra modules}
Let $ M $ be a $ \Lambda $-module of dimension vector $ \nu $.  By the definition of $\xi_M $ and the map $ D $, we have that
$$
D(\xi_M) =  \sum_{\vi} \chi(F_\vi(M))\, D_\vi
$$
In the previous section (Corollary \ref{co:DDH}), we saw that the measure $ D(b_Z) $ of an MV basis vector equals the asymptotics of sections of line bundles on $ Z$.  In a similar fashion, we will now explain that $ D(\xi_M) $ can also be regarded as an asymptotic.

Consider the algebra $ \Lambda[t] := \Lambda \otimes_\C \C[t] $.  We define
$$
\mathbb G_{\mu}(M[t]/t^n) = \bigl\{ N \subset M \otimes \C[t]/t^n : N \text{ is a $\Lambda[t]$-submodule, }\, \dimvec N = \mu \bigr\}
$$
We will record the information of the Euler characteristics of these varieties as an element of $ \Dist(\tR) $ by
$$ [H^\bullet(\mathbb G(M[t]/t^n))] = \sum_{\mu \in Q_+} \chi(\mathbb G_{\mu}(M[t]/t^n)) \; \delta_{\mu} $$

\begin{Theorem} \cite[Theorem 11.4 and Lemma 12.3]{mvbasis} \label{th:limLambda}
	For any $ \Lambda$-module $ M$, we have $$ D(\xi_M) = \lim_{n \rightarrow \infty} \frac{1}{n^{\dim M}} (\tau_n)_* [H^\bullet(\mathbb G(M[t]/t^n))] $$
\end{Theorem}

\subsection{A conjecture and symplectic duality}
Suppose that $ Y \in \Irr \Lambda(\nu) $, with general point $M$, and $ Z $ is a stable MV cycle, such that $ c_Y = b_Z $ (there are many such pairs, conjecturally).  Then $ D(c_Y) = D(b_Z) $.  Via Theorem \ref{th:limLambda} and Corollary \ref{co:DDH}, both sides are the asymptotics of $ T^\vee$-representations. So it is natural to expect equality before taking the limit.  If we further assume that the odd cohomology of $\mathbb G_\mu(M[t]/t^n) $ vanishes, this implies that there is an isomorphism of $ T^\vee$-representations,
\begin{equation} \label{eq:AllenConj}
	\Gamma(Z, \mathcal O(n)) \cong H^\bullet(\mathbb G(M[t]/t^n)), \quad \text{ for all } n \in \N,
\end{equation}
where $T^\vee $ acts on the right side through the decomposition $ \mathbb G(M[t]/t^n) = \sqcup \mathbb G_\mu(M[t]/t^n) $.


The left hand sides of (\ref{eq:AllenConj}) form the components of a graded algebra, so it is natural to search for a similar structure on the right hand sides.  After studying this question for some time, we are pessimistic about finding this algebra structure.  On the other hand, $ \C[Z] := \bigoplus_n \Gamma(Z, \mathscr L^{\otimes n}) $ is also a module over $ \C[\overline{S_+^\nu \cap S^0_-}] $. We believe that such a module structure naturally exists for the direct sums of the right hand side of (\ref{eq:AllenConj}). 

\begin{Conjecture}\label{conj:sheaf}
	\begin{enumerate}
		\item For any preprojective algebra module $ M $ of dimension vector $ \nu $, 
	$
	\bigoplus_{n \in \N}  H^\bullet(\mathbb G(M[t]/t^n))
	$
	carries the structure of a $C[\overline{S_+^\nu \cap S^{0}_-}]$-module.
	\item 
	When $b_Z = c_Y $ and $ M $ is a general point of $ Y$, then there is an isomorphism (\ref{eq:AllenConj}) of  $C[\overline{S_+^\nu \cap S^{0}_-}]$-modules.
	\end{enumerate}
\end{Conjecture}

This conjecture should be a manifestation of the symplectic duality between generalized affine Grassmannian slices and Nakajima quiver varieties.  In particular, Braverman-Finkelberg-Nakajima \cite{BFN2} proved that a generalized affine Grassmannian slice is the Coulomb branch associated to the quiver gauge theory defining the corresponding Nakajima quiver variety.  In \cite{BFNSpringer}, with Hilburn and Weekes, we developed a Springer theory for Coulomb branch algebras and proved a weak form of Conjecture \ref{conj:sheaf} for those modules $M $ which come from a representation of the undoubled quiver.

The symplectic singularity viewpoint is also a useful framework for thinking about our three bases.  In particular, following the philosophy of Braden--Licata--Proudfoot--Webster \cite{BLPW}, the MV cycles and quiver variety components can be categorified using category $\mathcal O $ for quantizations of affine Grassmannian slices and quiver varieties, respectively.  Moreover, these categories are closely related to categories of modules for KLR algebras \cite{YangianCatO}.  From this perspective, the failures of these bases to agree with the dual canonical basis in Theorem \ref{th:counterexample} can be attributed to the non-irreducibility of the character varieties of simple objects in these categories.

%
%
%
%
%
%
%

\section{Cluster structures}
\subsection{Cluster structures on $ \C[N]$}
Cluster algebras were defined by Fomin--Zelevinsky in order to understand the dual canonical basis of $ \C[N] $.  A cluster algebra is a commutative algebra $ A$ with a distinguished collection of ``clusters''.  Each \textbf{cluster} consists of an algebraically independent subset $ T = \{ x_1, \dots, x_n \} \subset A $, such that $ A \subset \C(x_1, \dots, x_n) $.  We pass from one cluster to another using an ``exchange procedure'' which removes one of the $ x_i$ and replaces it with a certain rational function.  A \textbf{cluster monomial} is a monomial in the variables in one cluster. 

Berenstein--Fomin--Zelevinsky \cite{BFZ} proved that $ \C[N] $ carries a cluster algebra structure.  Every reduced word $ s_{i_1} \cdots s_{i_m} = w_0 $ for the longest element of the Weyl group gives us a cluster $ T(\vi) $ (though these are not all the clusters).

Geiss--Leclerc--Schr\"oer \cite{GLSrigid} established the following beautiful result explaining how the theory of preprojective algebras provides an additive categorification of the cluster algebra structure on $ \C[N] $.

\begin{Theorem} \label{th:GLS1}  Assume that $ \fg $ is simply-laced.
	\begin{enumerate}
	\item A maximal rigid $ \Lambda$-module $ T $ gives a cluster with cluster variables $x_1 = \xi_{T_1}, \dots,$ $x_n = \xi_{T_n} $, where $ T_i $ are the distinct indecomposable summands of $ T $. 
	\item Every cluster is of this form and all cluster monomials lie in the dual semicanonical basis.
	\item The exchange relations in $ \C[N] $ comes from short exact sequences in $ \Lambda\mmod$.
	\end{enumerate}
\end{Theorem}

On the other hand, Kang--Kashiwara--Kim--Oh \cite{KKKO} proved that the categories of $R_\nu$-modules provide a monoidal categorification of the cluster algebra $\C[N]$.  In particular they proved the following result, which was obtained at around the same time by Qin \cite{Fan17}.

\begin{Theorem} \label{th:KKKO}
	Every cluster monomial in $\C[N]$ lies in the dual canonical basis.
\end{Theorem}

Together, Theorems \ref{th:GLS1} and \ref{th:KKKO} imply that the dual semicanonical and canonical bases contain many common elements, since they both contain all cluster monomials.

\subsection{$g$-vectors}

Fix a cluster $ T = \{x_1, \dots, x_n \} $ in a cluster algebra $ A$.  Let $ u \in A $ be a cluster monomial.  Fomin--Zelevinsky \cite{cluster4} defined combinatorially the $ g$-vector $ g_T(u) \in \Z^n $ of $ u $ using a mutation procedure.  In \cite{DWZ}, Derksen--Weyman--Zelevinsky proved that this $g$-vector encodes the ``leading monomial'' appearing in $u $.  In this case, we say $ u $ is $ g$-pointed.  As we vary the cluster $ T $, the data of these $ g_T(u) $ defines a tropical point in the Langlands dual cluster $ X $ variety, as studied by Fock-Goncharov \cite{FockG}.

 The following observation is due to Genz--Koshevoy--Schumann \cite[Section 6]{GKS}.
\begin{Proposition} \label{pr:gLusztig}
	Let $ \vi $ be a reduced word for $ w_0 $, giving a cluster $ T(\vi)$.  Let $ u $ be a cluster monomial.  Then $ g_{T(\vi)}(u) $ agrees with the $\vi$-Lusztig data of $ u$, up to a simple linear change of coordinates.  
\end{Proposition}
In this way, we see explicitly how the information of all these $ g$-vectors is the same information as the MV polytope (this is also closely related to Corollary \ref{co:GS}).

In the setting of the dual semicanonical basis, this can be generalized as follows.  Let  $ T $ be a maximal rigid $ \Lambda$-module and let $ M \in \Lambda\mmod$. Geiss--Leclerc--Schr\"oer \cite{GLSgeneric} defined $ g_T(M) \in \Z^n $ using homological algebra, extending the above notion of $g$-vector.  Moreover, they proved that $ \xi_M$ is $ g$-pointed in each cluster.

\subsection{Theta basis}
A cluster algebra that contains finitely many clusters is called \textbf{finite type}.  The cluster algebra $ \C[N] $ is of finite type only when $ G = \SL_2, \SL_3, \SL_4, \SL_5$.  When a cluster algebra $A$ is not finite type, then the cluster monomials do not span $ A$.  It was a longstanding open problem to extend the set of cluster monomials to a basis for $ A$.  This problem was solved by the following remarkable theorem of Gross--Hacking--Keel--Kontsevich.

\begin{Theorem} \label{th:GHKK}
	Let $ A $ be a cluster algebra, satisfying some hypotheses (which hold for $ \C[N]$).  There is a natural basis for $ A$, called the theta basis, extending the set of cluster monomials.  This theta basis is parametrized by the set of tropical points of the Langlands dual cluster $ X $ variety.  Moreover, each theta basis element is $g$-pointed in each cluster.
\end{Theorem}

Combining Theorem \ref{th:GHKK} with Proposition \ref{pr:gLusztig} or Corollary \ref{co:GS}, we get a natural parametrization of the theta basis of $ \C[N] $ by $ B(\infty) $.  In \S \ref{se:Binf}, we saw that biperfect bases are parametrized by $ B(\infty)$, so it is natural to ask the following.

\begin{Question}
	Is the theta basis for $\C[N] $ a biperfect basis?
\end{Question}

\subsection{Cluster structure on the MV basis}
From Theorems \ref{th:GLS1} and \ref{th:KKKO}  the dual semicanonical and dual canonical bases for $ \C[N] $ contain all cluster monomials.   In another direction, Qin \cite{Fan19} studied bases which are $g$-pointed in each cluster and gave a description of the set of all such bases.  In particular, he showed that all such bases contain all the cluster monomials.

This motivates the following conjecture.
\begin{Conjecture} \label{co:MVcluster}
	The MV basis for $ \C[N] $ contains all cluster monomials.  Moreover its elements are $g$-pointed in each cluster.
\end{Conjecture}

Note that the conjecture would imply that the MV basis and dual semicanonical basis agree for $\SL_5 $ which is not known.

As evidence for this conjecture, let us mention that the first counterexamples (in both $\SO_8$ and $ \SL_6$) to the equality of the MV and dual semicanonical bases (see Theorem \ref{th:counterexample}) occur for the square of the simplest basis element which is not a cluster monomial.

Baumann--Gaussent--Littelmann proved this conjecture for certain clusters.

\begin{Theorem} \cite[Prop 7.2]{BGL}
	If a reduced word $ \vi $ for $w_0$ satisfies a certain condition (which holds for all reduced words in small rank), then all cluster monomials in the cluster $ T(\vi) $  lie in the MV basis.
	\end{Theorem}

At the moment, we are far from Conjecture \ref{co:MVcluster}, but thinking about this conjecture motivates the following questions.
\begin{Question}
	\begin{enumerate}
		\item
	Is there a refinement of the notion of biperfect basis which would imply that such a basis contains all cluster monomials?
	\item
What do cluster exchange relations correspond to geometrically?  Which collections of MV cycles form clusters?
	\end{enumerate}
	\end{Question}


Finally, we close with the following wild conjecture.
\begin{Conjecture}
	The MV and theta bases for $ \C[N] $ coincide.
\end{Conjecture}
We have three pieces of weak evidence for this conjecture.  First, the way in which the MV, dual canonical, and dual semicanonical bases differ in $ \C[N] $ for $ \SL_6 $ in (\ref{eq:bcd}) is very reminiscent of the way in which the theta, dual canonical, and generic bases differ for rank 2, affine type cluster algebras.  Second, the construction of the MV and theta bases are both related to the  geometry of loop spaces.  Finally, the trichotomy of bases studied here seems to match the trichotomy of bases for cluster algebras, see \cite{Fan19}.

\subsection{Acknowledgements}
I would like to thank Pierre Baumann and Peter Tingley for many years of fruitful collaboration and friendship.  I also thank Alexander Braverman, Michel Brion, Elie Casbi, Anne Dranowski, Tom Dunlap, Michael Finkelberg, Stephane Gaussent, Justin Hilburn, Allen Knutson, Bernard Leclerc, Peter Littelmann, Calder Morton-Ferguson, Dinakar Muthiah, Ben Webster, Alex Weekes, and Harold Williams for collaboration and useful conversations.  Finally, I would like to thank Arun Ram his three mountains analogy.

This work was partially supported by an NSERC Discovery Grant.

\end{document}